\documentclass{article}
\usepackage{amsmath}
\usepackage{amsfonts}
\usepackage{amssymb}
\usepackage{graphicx}
\usepackage{booktabs}
\usepackage{cleveref}
\usepackage{tikz}

\makeindex

\newcommand{\proof}{\noindent {\textbf{Proof}}:\ }
\def\qed{{\hspace*{0mm}\hfill $\Box $\vspace{3mm}}}

\def\BB{\{0,1\}}

\def\RR{\mathbb R}

\def\bx{\mathbf{x}}
\def\bc{\mathbf{c}}

\def\by{\mathbf{y}}

\def\bZ{\mathbf{Z}}
\def\bf{\mathbf{f}}

\def\cG{\mathcal{G}}

\def\cL{\mathcal{L}}

\def\cA{\mathcal{A}}

\def\cE{\mathcal{E}}

\def\cQ{\mathcal{Q}}

\def\bP{\mathbf{P}}
\def\bQ{\mathbf{Q}}

\newcommand{\ga}{\alpha}
\newcommand{\gb}{\beta}

\newcommand{\gO}{\Omega}

\def\P{\mathcal{P}}

\def\P{\mathbb{P}}

\newtheorem{theorem}{Theorem}

\newtheorem{corollary}{Corollary}

\newtheorem{lemma}[theorem]{Lemma}

\newtheorem{remark}{Remark}

\begin{document}
	
	\title{Boole's probability bounding problem, linear programming aggregations, and nonnegative quadratic pseudo-Boolean functions}

	\author{
		Endre Boros\footnote{Email: endre.boros@rutgers.edu} and Joonhee Lee\footnote{Email: jlee15@pace.edu}\\
\\
		RUTCOR and MSIS Department\\
		Rutgers Business School\\
		Rutgers University\\
\\
		Management and Management Science Department\\
		Lubin School of Business\\
		Pace University\\
	}

	\date{January 8, 2025}
	
	\maketitle
	
	\begin{abstract}
		Hailperin (1965) introduced a linear programming formulation to a difficult family of problems, originally proposed by Boole (1854,1868). Hailperin's model is computationally still difficult and involves an exponential number of variables (in terms of a typical input size for Boole's problem). Numerous papers provided efficiently computable bounds for the minimum and maximum values of Hailperin's model
		by using aggregation that is a monotone linear mapping to a lower dimensional space. In many cases the image of the positive orthant is a subcone of the positive orthant in the lower dimensional space, and thus including some of the defining inequalities of this subcone can tighten up such an aggregation model, and lead to better bounds. Improving on some recent results, we propose a hierarchy of aggregations for Hailperin's model and a generic approach for the analysis of these aggregations.  We obtain complete polyhedral descriptions of the above mentioned subcones and obtain significant improvements in the quality of the bounds. 
	
	\end{abstract}

\section{Introduction}\label{intro}

In this paper we study a class of linear programming aggregations, motivated by a difficult family of problems introduced by Boole \cite{Boole1854,Boole1868} (see also \cite{Boole1952}). In the most studied ``union-bounding'' variant, we would like to know how large or small can be the probability of the union of $n$ events, given that we know the probabilities of these events and their pairwise intersections. Hailperin \cite{Hailperin65} introduced a linear programming formulation for this problem, involving $2^n$ variables. In his formulation the maximum and minimum values of the linear objective function are the best possible (tight) bounds for Boole's union-bounding problem. It turns out that these problems are computationally hard. For instance, the feasibility of Hailperin's model is NP-hard (see e.g., \cite{KavPap90}) and the column generation technique also leads to NP-hard optimization problems (see \cite{DezaLau97,HanJauNgu99}). Numerous papers studied these problems, and proposed computationally efficient bounds (that are of course not as tight as Hailperin's). Many of these techniques are using aggregations of Hailperin's model (though not all say this explicitely; for instance a large body of research considers the so called ``binomial moments'' problem, which in fact is an aggregation of Hailperin's model). We refer the reader to the introductions of e.g., \cite{Boros14,Prekopa88} for a more thorough overview of related research.

We view aggregations as monotone linear mappings $\cL:\RR_+^M\to \RR_+^N$, where typically $N\ll M$. We say that such a mapping is an aggregation of a linear programming problem with a polyhedron $\bP\subseteq \RR_+^M$ as feasible solutions and $\bc\in\RR^M$ as linear objective function, 
if we can construct a polyhedron $\bQ\subseteq \RR_+^N$ and a linear objective function $\bf\in \RR^N$ satisfying the following property:
\[
\forall~ \bx\in \bP \text{ we have } \cL(\bx)\in \bQ \text{ and } \bc^T\bx = \bf^T\cL(\bx).
\]
Furthermore we say that $\cL$ is a \emph{faithful aggregation} if in addition we have
\[
\forall~ \by\in \bQ ~~\exists~ \bx\in \bP ~\text{ such that }~ \by=\cL(\bx).
\]
It follows that an aggregation provides a relaxation for a linear program, while a faithful aggregation\footnote{Note that the LP optimum of a faithfull aggregation is a best possible bound for the original problem over all aggregations. Some authors use different terminology, including \emph{achivable} or \emph{optimal} for such bounds, see e.g., \cite{yang2017arxiv}.} is equivalent to it.

In this paper  we consider a hierarchy of aggregations of Hailperin's LP  formulation for bounding the union of $n$ events. The first trivial level $\cL_0$ with $N=n$ leads to the well known binomial moment problem. It is well known that this problem is an aggregation of Hailperin's model (see e.g., \cite{Prekopa88}), though the fact that $\cL_0(\RR_+^M)=\RR_+^n$ needs a (simple) proof, which we give in Section \ref{binomialmoments}. The  next level $\cL_1$ with $N=n^2$ coincides with the aggregation introduced by Pr\'ekopa and Gao \cite{PG}. For this case Qiu, Ahmed, and Dey \cite{qiu2016} noted that $\cL_1(\RR_+^M)\neq \RR_+^N$, and provided a complete polyhedral description for $\cL_1(\RR_+^M)$. In Section \ref{PGsection} we reprove their result with a new proof that allows us to introduce our technique connecting the valid inequalities of this cone to the cone of linear nonnegative \emph{pseudo-Boolean functions} (real valued functions in binary variables). The next level $\cL_2$ with $N=O(n^3)$ coincides with the aggregation proposed by Yang, Alajaji, and Takahara \cite{yang16,yang1}. While they noted that $\cL_2(\RR_+^M)\neq \RR_+^N$ and provided a number of valid inequalities for this subcone, they could not describe this subcone fully. Our main result is a complete polyhedral characterization of this subcone that we give in Section \ref{Yang-et-al}. In particular, we prove that the valid inequalities for this cone are in a one-to-one connection with the cone of quadratic nonnegative pseudo-Boolean functions. In fact, if we introduce all such valid inequalities, then we obtain a faithful aggregation of Hailperin's model in $N=O(n^3)$ dimension. Finally, using a polynomial subfamily of the described valid inequalities, we obtain polynomial time computable bounds for Boole's union-bounding problem, that are significantly better than other known efficiently computable bounds in most numerical examples we tested. 
	
\section{Notation and terminology}\label{Notation}

Let $n$ and $m$ denote positive integers and $V=[n]=\{1,2,...,n\}$ denote the set of indices. We denote by $\P$ a probability measure on a (possibly finite) probability space and use $\cA=\{A_i \mid i\in V\}$ to denote a collection of $n$ events in this space. We use $\cQ=\{Q\subseteq V \mid 0\leq |Q|\leq m\}$ to denote all subsets of indices of ``small'' size (we assume $m$ is a constant with respect to $n$), and use
$p_Q\in [0,1]$, $Q\in \cQ$ to denote probability values. Boole's union-bounding problem asks for the maximum and minimum values for $\P(A_1\cup\dots\cup A_n)$ subject to the constraints
\begin{equation}\label{e-PQ}
\P\left(\bigcap_{i\in Q} A_i\right) ~=~p_Q   ~~~\forall ~Q\in\cQ.
\end{equation}
Note that the intersection of an empty family of events is considered to be an event that is always true. In accordance, we define $p_\emptyset=1$.
Most papers consider a special case of these problems, corresponding to $m=2$. In the sequel we also focus on this special case. 

It is well known that $n$ events divides the probability space into $2^n$ \emph{atomic events}
\[
	\omega_S ~=~ \bigcap_{i\in S} A_i \cap \bigcap_{i\not\in S} \bar{A}_i  ~~~ S\subseteq V.
\]
Hailperin proposed to view
\[
	x_S ~=~ \P (\omega_S) ~~~\text{ for all }~~~ S\subseteq V	
\]
as unknown variables, and formulated the union bounding problem (in case of $m=2$) as the following linear optimization problem 
\begin{equation}\label{e-Hail0}
		\begin{array}{rrcll}
			\max/\min&\displaystyle\sum_{S \subseteq V\atop S \neq \emptyset}x_S &&&\\
			s.t. &&&&\\
			&\displaystyle\sum_{S\subseteq V\atop  S \supseteq Q} x_S &=& p_Q& ~~\forall~~ Q\subseteq V, 0\leq|Q|\leq 2,\\*[7mm]
			& x_S &\geq & 0 & ~~\forall~~ S\subseteq V.
		\end{array}
\end{equation}
Let us denote by $LB(H)$ and $UB(H)$ the optimum values of the minimization and maximization, respectively, of Hailperin's problem \eqref{e-Hail0}.
	
\bigskip
\noindent
Let us add here a few remarks, to enhance clarity of our following discussions:
\begin{itemize}
\item  Note that $LB(H)$ and $UB(H)$ are \emph{best possible} (also called sharp, optimal, or achievable bounds, see e.g., \cite{Boros14,yang2017arxiv}), meaning precisely that for an arbitrary real $Z$ with $LB(H)\leq Z\leq UB(H)$ there exists a probability space with events $A_i$, $i=1,...,n$ such that $Z=\P\left(\displaystyle\bigcup_{i=1}^n A_i\right)$ and equations \eqref{e-PQ} hold.
\item From a computational point of view, for a reasonably bounded rational input the input size is polynomially bounded by $n$, while the number of variables in \eqref{e-Hail0} is exponential in $n$. This in itself of course is not an indication of computational hardness. However, it was shown that trying to use column generation leads, in every iteration, to a binary unconstrained quadratic optimization, that in general is known to be NP-hard \cite{JauHanPog91}. Furthermore, the separation problem for the dual of the maximization in Hailperin's model is equivalent with a membership question for the so called cut-polytope, a well-known NP-complete problem \cite{DezaLau97}. It was also shown that the feasibility of \eqref{e-Hail0} belongs to the family of probabilistic satisfiablity problems, that are also known to be NP-hard, in general \cite{KavPap90}. These indications are not a definitive proof of computational hardness. In particular, it is easy to provide a feasible input for this problem, and then the complexity of the optimization problems \eqref{e-Hail0} is open. Nevertheless, the perceived hardness of \eqref{e-Hail0} motivated a large amount of research to provide bounds that are close to $LB(H)$ from below and to $UB(H)$ from above and are computable in polynomial time in terms of the input size. 
\item We can observe that variable $x_\emptyset$ appears only once, in the equality corresponding to $p_\emptyset=1$. Thus, we can eliminate both $x_\emptyset$ and the first equality from our formulation yielding a somewhat simplified formulation:	
\begin{equation}\label{e-Hail}
			\begin{array}{rrcll}
				\max/\min
				&\displaystyle\sum_{S \subseteq V\atop S \neq \emptyset}x_S &&&\\
				s.t. &&&&\\
				&\displaystyle\sum_{S\subseteq V\atop  S \supseteq Q} x_S &=& p_Q& ~~\forall~~ Q\subseteq V,  1\leq |Q|\leq 2,\\*[7mm]
				& x_S &\geq & 0 & ~~\forall~~ \emptyset\neq S\subseteq V.
			\end{array}
\end{equation}
This simplification may lead to a feasible formulation \eqref{e-Hail}, even when problem \eqref{e-Hail0} is infeasible \cite{Boros14}. Still, denoting by $Z^{\max}$ and $Z^{\min}$ the optimum values in \eqref{e-Hail}, we have the relations that 
\eqref{e-Hail0} is infeasible if and only if $Z^{\min} >1$, and for feasible inputs we have $LB(H)=Z^{\min}$ and $UB(H)=\min\{1,Z^{\max}\}$.
Thus, problems \eqref{e-Hail0} and \eqref{e-Hail} are equivalent. 
\end{itemize}
	
There is a very large body of research on providing good bounds for \eqref{e-Hail}, mostly by considering a relaxation that is either has polynomial dimension in $n$ or is at least polynomially computable, due to some structural property. We refer the reader to \cite{Boros14} for a concise summary of these results. 
Many of these results are based on relaxations of \eqref{e-Hail} obtained by aggregating some of the variables of the original formulation (and sometimes, aggregating also some of the equalities of \eqref{e-Hail}). 
	
	
In this paper we focus on aggregations of the variables of \eqref{e-Hail} of the following type:
	Consider a hypergraph $\cE=\{E_1,\ldots,E_N\}\subseteq 2^\gO$ where $\gO=2^V\setminus \{\emptyset\}$, and associate to it a new set if variables
	\begin{equation}\label{e-yk}
		y^k ~=~ \sum_{S\in E_k} x_S ~~~\text{ for } k=1,...,N
	\end{equation}
that are aggregations of the original variables of Hailperin's model.
	This in fact defines a monotone linear mapping $L_\cE: \RR^\gO \mapsto \RR^N$. 
	In each of the following special cases, we choose $\cE$ in such a way that the objective function and the equality constraints in \eqref{e-Hail} are ``easily representable'' (sometimes after applying additional aggregations of some of the equalities).  
	
	From the above definition it is also immediate that $L_\cE(\RR_+^\gO) \subseteq \RR_+^N$. A sometimes unrecognized property in the literature is that we may have $L_\cE(\RR_+^\gO) \neq \RR_+^N$. Thus, the obtained relaxation could frequently be tightened with additional inequalities corresponding to the facets of the cone $L_\cE(\RR_+^\gO)$. 
	
	
	\bigskip
	
	In the subsequent sections variables will be associated to various subfamilies of $\gO$. One of the dividing parameters used many times is the cardinality of the subsets. Accordingly, we introduce
	\[
	\gO^\ell ~=~ \{S\in \gO\mid |S|=\ell\}
	\]
	for $\ell=1,...,n$, and therefore we have the partition
	$\gO ~=~ \bigcup_{\ell=1}^n \gO^\ell$. 
	
	\bigskip

\section{The Binomial Moment Problem}\label{binomialmoments}
	
	Let us first consider the aggregation corresponding to the hypergraph $\cE^0=\{E^0_k\mid k=1,...,n\}$, where $E^0_k=\{S\subseteq V\mid |S|=k\}$ for $k=1,...,n$. This aggregation yields the well known \emph{binomial moment problem}:
	
	\begin{equation}\label{e-binomial}
	\max \left\{
	\sum_{k=1}^n y^k \left|
	\begin{array}{rcl} 
		\displaystyle\sum_{k=1}^n ky^k &=& B_1\\
		\displaystyle\sum_{k=1}^n \binom{k}{2} y^k &=& B_2\\
		y^k &\geq & 0 ~~\text{ for }~~ k=1,...,n.
	\end{array}
	\right.\right\}
	\end{equation}
	where $B_1=\sum_{ Q \in \gO^1} p_Q$ and $B_2=\sum_{Q \in \gO^2} p_Q$ are the first two binomial moments of the considered event system.

	Numerous papers dealt with this problem (\cite{BP1,CE52,DW1967,DeCaen99,Galambos,Hoppe1990,KAT,Kwerel75a,Kwerel75b,Prekopa88,SPS80,yang16}). It is well-known that for this case both the published upper and lower bounds are best possible, based on the aggregated information, see e.g. \cite{Prekopa88} for an LP analysis that implies this readily.
	In this aggregation the nonnegativity of the $x_S$, $S\in\gO$ variables are simply replaced by the nonnegativity of the aggregated $y^k$, $k=1,...,n$ variables. It turns out that in this case this is correct: this aggregation maps the positive orthant to the entire positive orthant in the lower dimension. While this can be derived by the approach of \cite{Prekopa88}, we include here a short proof for completeness.
	
	\begin{theorem}\label{t-binomial}
		We have 
		\[
		L_{\cE^0} \left(\RR_+^\gO\right) ~=~ \RR_+^n.
		\]
	\end{theorem}
	
	\proof
	Assume $\ga_0,...,\ga_1$ are real numbers such that 

	\begin{equation}\label{nonnegcone}
	\ga_0 +\sum_{k=1}^n \ga_ky^k ~\geq ~ 0
	\end{equation}
	whenever $\bx=(x_S\mid S\in\gO)\in \RR_+^\gO$. 
	
	Since $\bx=0\in \RR_+^\gO$ is a nonnegative vector, for which $\by=0$, the inequality $\ga_0\geq 0$ follows immediately. 
	
	Let us next consider vectors $\bx\in \RR_+^\gO$, in which $x_S=0$ for all $S\in \gO$, $|S|\neq k$. For these we have $y^j=0$ for all $j\neq k$, and $y^k$ can take arbitrary positive values. Since we must have $\ga_0+\ga_k\cdot y^k \geq 0$ for all these substitutions by the validity of \eqref{nonnegcone}, we must have $\ga_k\geq 0$, too. Because all coefficients in \eqref{nonnegcone} must be nonnegative, \eqref{nonnegcone} is a trivial consequence of the nonnegativity of the $y^k$, $k=1,...,n$ variables.
	\qed

		Let us add that \Cref{t-binomial} implies that for any other objective function in Hailperin's model that are linearly expressible in the aggregated $(y^k\mid k=1,...,n)$ space the LP optimums of \eqref{e-binomial} with this objective function are best possible (optimal) bounds with respect to the information represented by the binomial moments $B_1$ and $B_2$. For instance, for any threshold $1\leq t\leq n$ we can derive best possible bounds (with respect to the input information $B_1$ and $B_2$) for the probability that at least $t$ events happen (where $t=1$ corresponds to the union considered above) by using the objective function $\sum_{k=t}^ny^k$ in \eqref{e-binomial}.

	\bigskip

	\section{The Aggregation by Pr\'ekopa and Gao }\label{PGsection}
	
	Pr\'ekopa and Gao \cite{PG} introduced an aggregation of \eqref{e-Hail} that corresponds to the case where 
	\[
	\cE^1=\{E^k_i\mid i\in V,~ k=1,...,n\}, 
	\]
	and 
	\[
	E^k_i=\{S\in \gO^k \mid i\in S\}
	\]
	for all $1 \le k \le n$ and $ i \in V$.
	Thus, $L_{\cE^1}$ is a mapping from $\RR_+^\gO$ to the space of the $n^2$ variables defined by 
	\begin{equation}\label{PG-e}
		y^k_i ~=~ \sum_{S\in\gO^k\atop i\in S} x_S ~~~~\text{ for all } i\in V \text{ and } k=1,...,n.
	\end{equation}
	
	Using these definitions, the model they considered is the aggregation
	\begin{equation}\label{PG-model}
		\begin{array}{crcll}
			\max ~(\min) & \displaystyle\sum_{k=1}^n \frac{1}{k}\sum_{i\in V} y^k_i \\*[5mm]
			s.t. \\
			&\displaystyle\sum_{k=1}^n y^k_i &=& p_{\{i\}} &\forall ~~ i\in V,\\
			&\displaystyle\sum_{k=2}^n (k-1)\cdot y^k_i &=& \displaystyle\sum_{R\in\gO^2\atop i\in R} p_R&\forall~~ i\in V\\*[7mm]
			& y^k_i&\geq&0&
			\begin{array}{c}\forall~~ i\in V, \text{ and}\\ k=1,...,n.\end{array}
		\end{array}
	\end{equation}
	They proved that this maximization problem has the same maximum value as the upper bound by \cite{Kwerel75b} and the minimization problem has the same value as the lower bound by \cite{KAT}. Let us denote by $PG$ the linear programming problem formed by \eqref{PG-model} and by $LB(PG)$ and $UB(PG)$ the corresponding minimum and maximum values.
	
	We can observe that in this case
	$L_{\cE^1}(\RR_+^\gO)\neq \RR_+^{n^2}$, 
	thus formulation \eqref{PG-model} could potentially be tightened up (see also \cite{qiu2016} and \cite{yang16}).
	
	To see this, let us observe first that by \eqref{PG-e} we have the equalities
	\[
		\sum_{{S\subseteq V\atop |S|=k}\atop i\not\in S} x_S~~=~~ \sum_{j\in V\setminus\{i\}} y^k_j ~~-~~(k-1)\cdot y^k_i 
	\]
	for all $i\in V$ and $k=1,...,n$. 
	Thus, the inequalities
	\begin{equation}\label{Zi-bar}
	\sum_{j\in V\setminus\{i\}} y^k_j ~~-~~(k-1)\cdot y^k_i ~~\geq~~0
	\end{equation}
	are valid inequalities for $L_{\cE^1}(\RR_+^\gO)$, for all $i\in V$ and $1\leq k \leq n$, and they are not linear consequences of the nonnegativities in \eqref{PG-model}. 
	
	Let us introduce the convex cone
	\[
	U^1 ~~~=~~~ 
	\left\{
	\by\in \RR_+^{n^2} \left|
	\begin{array}{rcl}
		y^k_i &\geq& 0 \\*[5mm]
		\displaystyle\sum_{j\in V\setminus\{i\}} y^k_j - (k-1)\cdot y^k_i &\geq& 0 
	\end{array} 
	~~~
	\begin{array}{c}
		\forall~ i\in V \text{ and}\\ k=1,...,n.\\ $~~~$
	\end{array} 
	\right.\right\}
	\]
	Then we have the relations 
	\begin{equation}\label{e-L-U}
		L_{\cE^1}(\RR_+^\gO) ~\subseteq~ U^1 ~\subsetneq~ \RR_+^{n^2}.
	\end{equation}
	
	Let us pause for a moment, and note that by its definition $L_{\cE^1}(\RR_+^\gO)$ is a convex cone, and since it is the image by a linear mapping of a polyhedral region, itself is polyhedral. Our first result in this section is to provide a full characterization of its facets. Note that this result was published in \cite{qiu2016} with a different proof. We keep our proof, since it introduces some of the proof ideas that will be instrumental in the next level of aggregations considered in the next section. 
	
	For real numbers $\ga_0$, $\ga^k_i$, $i\in V$, $k=1,...,n$ we call the inequality 
	\begin{equation}\label{e-valid}
		\ga_0 ~+~ \sum_{k=1}^n\sum_{i\in V} \ga^k_i\cdot y^k_i ~~\geq~~ 0
	\end{equation}
	a \emph{valid inequality for $L_{\cE^1}(\RR_+^\gO)$} if it holds for all $\bx\in\RR_+^\gO$ under the substitutions \eqref{PG-e}.
	Our first observation is that facets of $L_{\cE^1}(\RR_+^\gO)$ (tightest valid inequalities that are not convex combinations of other valid inequalities) have a special format.
	
	\begin{lemma}\label{l-tightest1a}
		If \eqref{e-valid} is a facet of $L_{\cE^1}(\RR_+^\gO)$, then we must have $\ga_0=0$, and there exists an integer $1\leq k \leq n$ such that $\ga^\ell_i=0$ for all $i\in V$ and $\ell\neq k$. 
	\end{lemma}
	
	\proof
	Let us note first that for $\bx=0\in \RR_+^\gO$ we have $y^k_i=0$ for all $i\in V$ and $k=1,...,n$. Thus, $\ga_0\geq 0$ follows. 
	
	We claim next, that for an arbitrary value of $k$ ($1\leq k\leq n$) the inequality
	\[
	\sum_{i\in V} \ga^k_i\cdot y^k_i ~~\geq~~0
	\]
	is also a valid inequality for $L_{\cE^1}(\RR_+^\gO)$. For this let us denote the left hand side of the above inequality as $F^k(\bx)$ for $\bx\in \RR_+^\gO$, $k=1,...,n$, assuming the substitutions \eqref{PG-e}. Let us also note that this expression depends only on  variables $x_S$ with $|S|=k$. Consequently we can assume that $x_S=0$ for all $S\in \gO$ with $|S|\neq k$. Assume then for a contradiction that the above inequality is not valid, that is that there exits a vector $\bx\in \RR_+^\gO$ with $x_S=0$ for all $S\in \gO$ with $|S|\neq k$ such that $F^k(\bx)<0$. For such a vector, and for any positive real $\lambda>0$ we have that the left hand side of \eqref{e-valid} is equal to
	\[
	\ga_0 + F^k(\lambda\cdot \bx) ~=~ \ga_0 + \lambda \cdot F^k(\bx)
	\]
	which for large $\lambda$ values would turn negative, contradicting the validity of \eqref{e-valid}.
	This contradiction, when applied for all values of $k$, and the nonnegativity of $\ga_0$ implies that \eqref{e-valid} is a linear consequence of the inequalities
	\[
	\sum_{i\in V} \ga^k_i\cdot y^k_i ~~\geq~~0 ~~~\forall k=1,...,n.
	\]
	\qed
	
	\begin{lemma}\label{l-tightest1b}
		Given reals $\ga_i$, $i\in V$, and an integer $1\leq k \leq n$, the inequality
		\begin{equation}\label{e-all-valid}
			\sum_{i\in V} \ga_i\cdot y^k_i ~~\geq~~ 0
		\end{equation}
		is valid for $L_{\cE^1}(\RR_+^\gO)$ if and only if for all subsets $S\in\gO^k$ we have 
		\[
		\sum_{i\in S} \ga_i ~~\geq~~ 0.
		\]
	\end{lemma}
	\proof
	By Lemma \ref{l-tightest1a} all facets of $L_{\cE^1}(\RR_+^\gO)$ have the form like in \eqref{e-all-valid}. Viewing the left hand side as a linear combination of the $x_S$, $S\in\gO^k$ variables after substitutions in \eqref{PG-e}, it is easy to see that 
	\[
	\sum_{i\in S} \ga_i
	\] 
	is the coefficient of $x_S$ in this linear combination, for all $S\in \gO^k$. Thus, if this quantity was negative for a subset $S\in \gO^k$, then changing the value of $x_S$ to a larger value, while not changing the values of the other $x_R$, $R\in \gO^k$, $R\neq S$ variables would strictly decrease (at a linear rate) the left hand side of \eqref{e-all-valid}. Since $x_S$ could take an arbitrarily large value, the claim follows. 
	\qed
	
	\bigskip
	
	The main result in this section (see also \cite{qiu2016}) is the following equality:
	
	\begin{theorem}\label{t-PG}
		\[
		L_{\cE^1}(\RR_+^\gO) ~=~ U^1.
		\]
	\end{theorem}
	
	\proof
	According to \eqref{e-L-U}, $U^1$ is not smaller than $L_{\cE^1}(\RR_+^\gO)$. Assume for a contradiction that $U^1\neq L_{\cE^1}(\RR_+^\gO)$, and choose a feasible solution $\by\in U^1\setminus 
	L_{\cE^1}(\RR_+^\gO)$. This implies that $\by$ must violate one of the facets of $L_{\cE^1}(\RR_+^\gO)$. According to Lemmas \ref{l-tightest1a} and \ref{l-tightest1b}, such a (violated) inequality should look like
	\begin{equation}\label{e-violated}
		\sum_{i\in V} \ga_i\cdot y^k_i ~~ < ~~ 0,
	\end{equation}
	for some positive integer $1 \le k \le n$.
	According to Lemma \ref{l-tightest1b}, permuting the coefficients in this inequality also yields a valid inequality for $L_{\cE^1}(\RR_+^\gO)$. Thus we can assume a permutation that yields the lowest left hand side in \eqref{e-violated}, i.e., we can assume that 
	\begin{equation}\label{sorted}
		\begin{array}{c}
			y^k_{1} ~\geq~ y^k_{2} ~\geq~ \cdots ~\geq~ y^k_{n}\\
			\ga_{1} ~\leq~ \ga_{2} ~\leq~ \cdots ~\leq~ \ga_{n}.
		\end{array}
	\end{equation}
	Inequality \eqref{e-violated} implies that we must have $\ga_{1} <0$. Furthermore by Lemma \ref{l-tightest1b} we get $\ga_{1}+ \cdots +\ga_{k}\geq 0$, implying $\ga_{k}>0$. 
	
	Note that after the sorting of the components as in \eqref{sorted}, if we decrease $\ga_{i}$ by epsilon, and increase $\ga_{j}$, $j>i$ by the same $\epsilon$, then the left hand side of \eqref{e-violated} can only decrease. Let us perform this operations, decreasing $\ga_1$ and increasing $\ga_j$ for $2\leq j< k$ until we arrive to the vector
	\[
	\begin{array}{rcl}
		\gb_{i} ~=~ 
		\begin{cases}
			\ga_{1} -(k-2)\ga_{k} +\ga_{2} + \cdots +\ga_{k-1} & \text{ if } i=1,\\
			\ga_{k} & \text{ if } 2\leq i \leq k,\\
			\ga_{i} & \text { if } i > k.
		\end{cases}
	\end{array}
	\]
	Observe finally that we have $\gb_{1}\geq -(k-1)\ga_{k}$, and $\gb_{i}\geq \ga_{k}$ for all $i=2,...,n$. Thus we can write the chain of inequalities
	\[
	0 > \sum_{i\in V}\ga_i\cdot y^k_i \geq \sum_{i\in V} \gb_i \cdot y^k_i \geq \ga_{k} \left(-(k-1)y^k_{1}+\sum_{i=2}^n y^k_{i} \right) ~\geq~ 0
	\]
	Here the last inequality follows since $\ga_{k}>0$ and $\by\in U^1$.
	The derived contradiction proves that $U^1=L_{\cE^1}(\RR_+^\gO)$, as claimed.
	\qed
	
Using Theorem \ref{t-PG} we can reformulate Pr\'ekopa-Gao model as follows:
	
	\begin{equation}\label{IPG-model}
		\begin{array}{crcll}
			\max ~(\min) & \displaystyle\sum_{k=1}^n \frac{1}{k}\sum_{i\in V} y^k_i \\*[5mm]
			s.t. \\
			&\displaystyle\sum_{k=1}^n y^k_i &=& p_{\{i\}} &~~\forall ~~ i\in V,\\
			&\displaystyle\sum_{k=2}^n (k-1)\cdot y^k_i &=& \displaystyle\sum_{R\in\gO^2\atop i\in R} p_R&~~\forall~~ i\in V\\*[7mm]
			&\displaystyle\sum_{j\in V\setminus\{i\}} y^k_j - (k-1)\cdot y^k_i &\geq& 0 &
			\begin{array}{c}\forall~~ i\in V, \text{ and}\\ k=1,...,n.\end{array}\\*[7mm]
			& y^k_i&\geq&0&
			\begin{array}{c}\forall~~ i\in V, \text{ and}\\ k=1,...,n.\end{array}
		\end{array}
	\end{equation}
	
	Let us denote by $IPG$ the linear programming problem formed by \eqref{IPG-model} and by $LB(IPG)$ and $UB(IPG)$ the corresponding minimum and maximum values.
Note that Yang et al. \cite{yang16} and Qiu et al. \cite{qiu2016} also considered $U^1$ and tightening the bound. They analyzed when the bound is less than or equal to 1 while our result holds for any feasible solution of \eqref{IPG-model}.

		\begin{table} [h]
	\small
			\centering
			\begin{tabular}{c c c}
				\toprule
				
				n & $\displaystyle \left( \frac{UB(PG)-UB(IPG)}{UB(PG)-UB(H)}\right) \times 100$ & $\displaystyle \left( \frac{LB(IPG)-LB(PG)}{LB(H)-LB(PG)} \right) \times 100$\\
				\midrule
	8	&	32.635\%	&	0.000\%	\\
	9	&	27.736\%	&	0.000\%	\\
	10	&	24.327\%	&	0.000\%	\\
	11	&	24.229\%	&	0.000\%	\\
	12	&	22.821\%	&	2.000\%	\\
	13	&	21.705\%	&	6.000\%	\\
	14	&	20.783\%	&	14.000\%	\\
	15	&	18.121\%	&	22.000\%	\\
	16	&	17.695\%	&	38.000\%	\\
	17	&	18.940\%	&	55.434\%	\\
	18	&	17.157\%	&	62.159\%	\\
	19	&	15.720\%	&	54.316\%	\\
	20	&	14.995\%	&	38.999\%	\\
	21	&	13.742\%	&	14.040\%	\\
	
				\bottomrule
			\end{tabular}
			\caption{Comparision of the PG and IPG bounds: We show the mean of the percentage relative error improvements over 50 randomly generated problems for each size. We write 0\% if both bounds are best possible. We did not test examples with $n>21$, because computing the $UB(H)$ and $LB(H)$ bounds are computationally intensive. 
				\label{tableIPGB}}
		\end{table}
	
	In Table \ref{tableIPGB} we compare these bounds. 
	In many examples, the addition of $U^1$ improves both the lower and upper bounds uniformly, sometimes quite significantly. Note that our computations are in line with the experience of \cite{qiu2016} that also showed that adding the inequalities of $U^1$ significantly improves on the PG bounds. 
	
	\bigskip
	
	\section{The Aggregation by Yang et al.}\label{Yang-et-al}
	
	In \cite{yang1} a novel aggregation model was considered that in our terms and notation corresponds to the case of 
	\[
	\cE^2 ~=~ \{E^k_Q \mid Q\in\gO^1\cup\gO^2,~~ |Q|\leq k\leq n\},
	\]
	where 
	\[
	E^k_Q ~=~ \{ S\subseteq V\mid Q\subseteq S, ~~ |S|=k\} ~~~\forall~ Q\in\gO^1\cup\gO^2,~~ |Q|\leq k\leq n \}.
	\]
	In other words, they consider the aggregation,  defined by the following mapping $L_{\cE^2}(\bx)=\by=(y^k_Q\mid Q\in\gO^1\cup\gO^2,~~ |Q|\leq k\leq n)$, where
	\begin{equation}\label{aijk}
		y^k_Q ~=~ \sum_{S\in\gO^k\atop Q\subseteq S} x_S ~~~\forall~ Q\in\gO^1\cup\gO^2,~~ |Q|\leq k\leq n.
	\end{equation}
	
	Note also that this notation is somewhat redundant, since we have the equalities
	\begin{equation}\label{yki}
		y^k_{\{i\}}~=~ \frac{1}{k-1} \sum_{Q\in \gO^2\atop i\in Q} y^k_Q
	\end{equation}
	for all $i\in V$ and $2\leq k\leq n$. Note finally that in general there are no other  linear relations between these variables.
	Thus in fact, this aggregation is a linear mapping $L_{\cE^2}$ into the space of dimension $N=n+(n-1)\cdot\binom{n}{2}$. 
	
	It is easy to verify that with this aggregation the objective function and the left hand sides of the equality constraints in \eqref{e-Hail} can be represented as linear expressions without any further aggregations of those equalities.
	
	\begin{align*}
		\displaystyle \sum_{S\in\gO} x_S ~=~& \displaystyle 
		\sum_{Q\in\gO^1} y^1_Q ~+~
		\sum_{k=2}^n \frac{1}{\binom{k}{2}}\sum_{Q\in\gO^2} y^k_Q&\rightarrow&~\max~ &&\label{obj}\tag{$E1$}\\
		\displaystyle \sum_{S\in\gO\atop i\in S} x_S ~=~& \displaystyle y^1_{\{i\}} ~+~ \sum_{k=2}^n \frac{1}{k-1} \sum_{Q\in\gO^2\atop i\in Q} y^k_Q &~=~&p_{\{i\}}&& \forall ~ i\in V\label{pi}\tag{$E2$}\\
		\displaystyle \sum_{S\in\gO\atop Q\subseteq S} x_S ~=~& \displaystyle \sum_{k=2}^n y^k_Q &~=~&p_Q&& \forall ~ 
		Q\in\gO^2 \label{pij}\tag{$E3$}\\
		&y^k_{\{i\}}&~\geq~&0&&\forall ~ \begin{array}{c}i\in V\\1\leq k\leq n\end{array} \label{Zi}\tag{$E4$}\\
		&y^k_Q&~\geq~&0&&\forall ~ \begin{array}{c}Q\in\gO^2\\2\leq k\leq n\end{array} \label{ZiZj}\tag{$E5$}
	\end{align*}
	leading to a quite tight polynomial size aggregation of our exponential sized original problem. Note that inequalities \eqref{Zi} include redundant ones, according to \eqref{yki}, but we keep this redundancy, since it will help us to realize a more general structure to these inequalities. 
	
	The nonnegativity of the $y^1_{\{i\}}$, $i\in V$ and $y^k_Q$, $Q\in\gO^2$ variables are implied by their definition \eqref{aijk}, since all $x_S$ variables are nonnegative. The authors of \cite{yang1} however noticed that in this case $L_{\cE^2}(\RR_+^\gO) \neq \RR_+^N$. In fact, 
	several other nonnegative combination of the $x_S$ variables are linear functions of $\by$. For instance, we can observe that for $i\in V$ and integers $2\leq k \leq n$ we have
	\[
	\sum_{S\in \gO^k\atop i\not\in S} x_S 
	~~=~~
	\frac{1}{\binom{k}{2}} 
	\sum_{Q\in\gO^2} y^k_Q
	- \frac{1}{k-1}
	\sum_{Q\in\gO^2\atop i\in Q} y^k_Q 
	~\geq~ 0
	\]
	and the nonnegativity of this expression is not a consequence of the nonnegativity of $\by$. Thus the nonnegativity of this expression could be added to the aggregation, tightening it up.
	In fact in \cite{yang1} the authors came up with the following additional nonnegative combinations of the $\bx$ variables that can be expressed as a linear function of the $\by$ variables when we use substitutions \eqref{aijk}.
	
	\[
	\sum_{S\in\gO^k \atop R\cap S =\emptyset} x_S 
	~=~ 
	\frac{1}{\binom{k}{2}} \left( \sum_{Q\in\gO^2} y^k_Q  - \frac{k}{2}\sum_{Q\in\gO^2\atop |R\cap Q|=1} y^k_Q  
	+  \frac{k(k-3)}{2}y^k_R \right)  
	~\geq~ 0
	\label{ZibarZjbar}\tag{$E6$}
	\]
	for all $R\in \gO^2$, and $2\leq k\leq n$. 
	
	\bigskip
	
	\[
	\sum_{{S\in \gO^k}\atop R\cap S=\{i\}} x_S
	~=~  
	\frac{1}{k-1} 
	\left(
	\sum_{Q\in\gO^2\atop  Q\cap R=\{i\}} y^k_Q - (k-2)y^k_R 
	\right) 
	~\geq~ 0
	\label{ZiZjbar}\tag{$E7$}
	\]
	for all $i\in V$, $i\in R\in \gO^2$, and $2\leq k\leq n$. 
	\bigskip
	
	\[
	\sum_{S\in\gO^k \atop |R\cap S|\in\{0,3\}}x_S ~=~ \frac{1}{\binom{k}{2}}\sum_{Q\in\gO^2} y^k_Q
	- \frac{1}{k-1} \sum_{Q\in\gO^2\atop |R\cap Q|=1} y^k_Q 
	+ \frac{k-3}{k-1}\sum_{Q\in\gO^2\atop Q\subseteq R} y^k_Q
	~\geq~ 0
	\label{ZiZjZl+ZibarZjbarZlbar}\tag{$E8$} 
	\]
	for all $R\in\gO^3$ and $2\leq k\leq n$.
	
	\bigskip
	
	\[
	\sum_{{{S\in\gO^k\atop R\cap S=\{i\}}\atop \textbf{or}}\atop R\setminus S=\{i\} }x_S 
	~=~
	\frac{1}{k-1}
	\sum_{Q\in\gO^2\atop Q\cap R=\{i\}} y^k_Q 
	+ y^k_{R\setminus\{i\}}
	-\frac{k-2}{k-1}\sum_{Q\in\gO^2\atop i\in Q\subseteq R} y^k_Q 
	~\geq~ 0
	\label{ZiZjbarZlbar+ZibarZjZl}\tag{$E9$}
	\]
	for all $i\in V$, $i\in R\in\gO^3$, and $2\leq k\leq n$. 
	
	\bigskip
	
	Let us remark that the polynomially sized linear programming aggregation \eqref{obj} - \eqref{ZiZjbarZlbar+ZibarZjZl} proposed by \cite{yang1} provides a bound for our original problem that is not weaker than the one proposed in the previous section. This is because problem \eqref{PG-model} and its strenghtening by $U^1$ are both aggregations of problem \eqref{obj} - \eqref{ZiZjbarZlbar+ZibarZjZl}. To see this, note that $y^1_i$, $i\in V$ have identical definitions in both aggregations (see \eqref{PG-e} and \eqref{aijk}), and for $y^k_i$, $i\in V$ and $k=2,...,n$ we have the equalities \eqref{yki}. Furthermore the inequalities defining $U^1$ are consequences of \eqref{Zi} - \eqref{ZiZjbar}.
	
	\bigskip
	
	\textbf{The question arise: do the inequalities \eqref{Zi}-\eqref{ZiZjbarZlbar+ZibarZjZl} describe the cone $L_{\cE^2}(\RR^\gO_+)$? } The answer is no. In what follows we provide an exact polyhedral description of $L_{\cE^2}(\RR^\gO_+)$.

	Let us start observing that the analogue of Lemma \ref{l-tightest1a} can be shown for this case, too.

	For real numbers $\ga_0$, $\ga^k_Q$, $Q \in \gO^1\cup\gO^2$, $k=1,...,n$ we call the inequality 
	\begin{equation}\label{e-valid-YAT}
		\ga_0 ~+~ \sum_{k=1}^n\sum_{Q\in\gO^1\cup\gO^2} \ga^k_Q\cdot y^k_Q ~~\geq~~ 0
	\end{equation}
	a \emph{valid inequality for $L_{\cE^2}(\RR_+^{\gO})$} if it holds for all $\bx\in\RR_+^{\gO}$ under the substitutions \eqref{aijk}.
	
	\begin{lemma}\label{l-tightest-2}
		If \eqref{e-valid-YAT} is a facet of $L_{\cE^2}(\RR_+^{\gO})$, then we must have $\ga_0=0$, and there exists an integer $1\leq k \leq n$ such that $\ga^\ell_Q=0$ for all $Q \in \gO^1\cup\gO^2$ and $\ell\neq k$. 
	\end{lemma}
	
	\proof
	The same poof idea we used for Lemma \ref{l-tightest1a} is working in this case, too. We omit to include this proof for brevity.
	
	\qed
	
	\bigskip
	
	Thus we can focus on subsets of a certain size. Let us also note that for sets $S\in\gO^1$ we have $x_S=y^1_S\geq 0$ by \eqref{Zi}, thus we can focus on set sizes $2\leq k\leq n$. Note also that all inequalities \eqref{Zi}-\eqref{ZiZjbarZlbar+ZibarZjZl} involve $y^k_Q$ variables for a fixed value of $k$. 
	
	For notational simplicity let us introduce $\gO^0=\{\emptyset\}$, set $\Sigma=\gO^0\cup\gO^1\cup\gO^2$, and define for $2\leq k\leq n$
	
	\begin{equation}\label{yk0}
		y^k_{\emptyset} ~=~ \sum_{S\in\gO^k} x_S ~=~ \frac{1}{\binom{k}{2}}\sum_{Q\in\gO^2} y^k_Q,
	\end{equation}
	where the second equality follows by the definitions of the $\by^k$ variables, see \eqref{aijk}.
	
	Let us now fix a $2\leq k\leq n$ value, and introduce $\bx^k=(x_S\mid S\in\gO^k)$ and $\by^k=(y^k_Q\mid Q\in \Sigma)$. Recall that the $y^k_Q$ variables are in fact (linear) functions of $\bx^k$ according to \eqref{aijk}, \eqref{yki}, and \eqref{yk0}. Thus, for a real vector $\ga\in\RR^{\Sigma}$ we say that the inequality
	\begin{equation}\label{valid}
		\sum_{Q\in\Sigma} \ga_Q\cdot y^k_Q ~\geq~ 0
	\end{equation}
	is \emph{valid} if it holds for all nonnegative $\bx^k\in \RR^{\gO^k}_+$.
	We denote the left hand side of inequality \eqref{valid} by $L_\ga(*)$, where $*$ indicates the set of variables we want to use. Thus inequality \eqref{valid}, as it is written, is the expression $L_\ga(\by^k)\geq 0$. We can also view the left hand side as a function of $\bx^k$, and then we write $L_\ga(\bx^k)$. Due to the relations \eqref{aijk},\eqref{yki}, and \eqref{yk0}, we can also write $L_\ga(y^k_Q\mid Q\in \gO^2)$. 
	Note finally that the same $\ga$ vector can be used as coefficients for different $k$ values, e.g., $L_\ga(\by^\ell)~=~L_\ga(\bx^\ell)$ for $\ell\neq k$ are also well defined.
	
	\bigskip
	
	In what follows we characterize all valid inequalities (for all $2\leq k \leq n$), and hence provide a complete description for $L_{\cE^2}(\RR^\gO_+)$. To arrive to such a complete description, we consider a mapping into the space of nonnegative quadratic pseudo-Boolean functions. 
	
	Introduce $\bZ=(Z_i\mid i\in V)\in \BB^V$, and associate to inequality \eqref{valid} the quadratic pseudo-Boolean function (or QPBF in short)
	\begin{equation}\label{F(Z)}
		F_\ga(\bZ) ~=~ \sum_{Q\in\Sigma} \ga_Q\cdot \prod_{i\in Q} Z_i,
	\end{equation}
	where we have $\prod_{i\in \emptyset}Z_i ~=~1$ by definition. 
	For a subset $S\subseteq V$ we denote by $\chi_S\in\BB^V$ its characteristic vector. 
	
	\begin{lemma}\label{qpbf-xs}
		Given a real vector $\ga\in\RR^\Sigma$ and a subset $\emptyset\neq S\subseteq V$, the real value $F_\ga(\chi_S)$ is the coefficient of $x_S$ in $L_\ga(\bx^{|S|})$. 
	\end{lemma}
	\proof
	Observe that $x_S$ appears in $y^k_Q$ with coefficient $1$ exactly when $Q\subseteq S$, and we have 
	\[
	F_\ga(\chi_S)~=~\sum_{Q\in\Sigma\atop Q\subseteq S} \ga_Q
	\]
	by \eqref{F(Z)}.
	\qed
	
	\begin{corollary}\label{c1}
		If $F_\ga(\bZ)\geq 0$ for all $\bZ\in\BB^V$, then $L_\ga(\by^k)\geq 0$ is a valid inequality for all $k=2,...,n$. 
	\end{corollary}
	\proof
	By Lemma \ref{qpbf-xs} the coefficients of the $x_S$, $S\in\gO$ variables are all nonnegative, since $F_\ga$ is a nonnegative function. Thus $L_\ga(\by^k)$ is a nonnegative real for all $\bx^k\in \RR^\Sigma_+$. 
	\qed
	
	Since linear combinations of QPBF-s is again a QPBF, the set of nonnegative QPBF-s form a convex cone (in the space of their coefficients, i.e., in dimension $1+n+\binom{n}{2}$). It is well-known that this cone is polyhedral, and its extremal rays are in a one-to-one correspondence with the facets of the cut polytope \cite{DeSimone89}. Let us denote by $\cG$ the finite set of extremal nonnegative QPBF-s, and thus $cone(\cG)$ is the set of nonnegative QPBF-s. 
	
	Let us now consider the polyhedral cone defined by the valid inequalities corresponding to nonnegative QPBFs. 
	
	\begin{equation}\label{e-U2}
		U^2  ~=~ 
		\left\{ \by ~\left|~
		\begin{array}{rlc}
			y^1_{\{i\}} &\geq~ 0 &\forall~ i\in V,\\*[5mm]
			\displaystyle\sum_{Q\in\Sigma} \ga_Q\cdot y^k_Q &\geq~ 0  &\begin{array}{c}
				\forall~\ga\in\RR^\Sigma \text{ and } 2\leq k\leq n\\
				\text{such that } F_\ga(\bZ)\in \cG\end{array} 
		\end{array}
		\right.\right\}
	\end{equation}
	
	Since by Corollary \ref{c1} all these inequalities are valid for $L_{\cE^2}(\RR^\gO_+)$, we have the relation
	\[
	L_{\cE^2}(\RR^\gO_+) ~\subseteq~ U^2
	\]
	
	In what follows we show that in fact we have equality here. The main difficulty in proving this claim stems from the fact that
	an inequality of the form \eqref{valid} may correspond to multiple QPBF-s $F_\ga(\bZ)$. This is because of the linear dependencies, like equations \eqref{yki} and \eqref{yk0}, we have in the $\by^k$ space. In fact, to a given QPBF we have several corresponding inequalities of the form \eqref{valid}, with different $k$ values. Some may be valid for $L_{\cE^2}(\RR^\gO_+)$ and some may not. More precisely, a real vector $\ga$ as coefficients in such an inequality may yield a valid one with variables $\by^k$, while if we replace $\by^k$ with $\by^\ell$ for some $\ell\neq k$, then the same expression may not be a valid inequality. Thus we need a mechanism that can transform a QPBF into another one, without changing a corresponding inequality (for a particular value of $k$). 
	
	A very helpful observation is the following claim:
	
	\begin{lemma}\label{identityk2}
		Given an integer $2\leq k\leq n$, let $\ga\in \RR^\Sigma$ be the coefficient vector of the QPBF
		\[
		\left(k-\sum_{i\in V}Z_i\right)^2 ~=~ k^2 -(2k-1)\cdot\sum_{i\in V} Z_i + 2\cdot\sum_{1\leq i<j\leq n} Z_i Z_j. 
		\]
		Then $L_\ga(y^k_Q\mid Q\in\gO^2)$ is the identically zero function.  
	\end{lemma}
	\proof
	Let us observe first that the equality in the above formula holds, because for binary variables we have $Z_i^2=Z_i$ for all $i\in V$.
	Thus we have
	\[
	L_\ga(\by^k) ~=~ k^2y^k_\emptyset -(2k-1)\sum_{Q\in \gO^1} y^k_Q + 2\sum_{Q\in \gO^2}y^k_Q.
	\]
	The claim now follows by elementary algebra, using the equalities \eqref{yki} and \eqref{yk0}.
	
	\begin{align*}
		L_\ga(y^k_Q\mid Q\in \gO^2) &=~ k^2\left(\frac{1}{\binom{k}{2}}\sum_{Q\in\gO^2} y^k_Q\right)- (2k-1)\sum_{i\in V}\left(\frac{1}{k-1}\sum_{Q\in\gO^2\atop i\in Q} y^k_Q\right)\\
		&~~~+2\sum_{Q\in\gO^2} y^k_Q
	\end{align*}
	Thus, for an arbitrary $Q\in\gO^2$ the coefficient of $y^k_Q$ in the above expression is
	\[
	k^2\frac{1}{\binom{k}{2}}-(2k-1)\frac{2}{k-1}+2~~~=~~~0.
	\]
	
	\qed
	
	Note that $L_\ga(y^\ell_Q\mid Q\in\gO^2)$ may not be identically zero if $\ell\neq k$.
	
	\bigskip
	
	Now, we are ready to prove our main result in this section.
	
	\begin{theorem}\label{main2}
		\[
		L_{\cE^2}(\RR^\gO_+) ~=~ U^2
		\]
	\end{theorem}
	
	\proof
	For $k=1$ we have $y^1_{i} = x_{\{i\}}$ for all $i\in V$, and clearly no other inequality of these variables can be minimally valid. 
	For $k\geq 2$, the containment``$\subseteq$" is implied by Corollary \ref{c1}. To see the reverse containment, let us assume that $\ga^*\in\RR^\Sigma$ is a real vector for which the inequality
	\[
	\sum_{Q\in \Sigma} \ga^*_Qy^k_Q ~\geq~ 0
	\]
	is valid, and for which $F_{\ga^*}(\bZ)$ is \textbf{not} a nonnegative QPBF. 
	Let us introduce 
	\[
	\gamma ~=~ \min_{\bZ\in\BB^V} F_{\ga^*}(\bZ) ~<~ 0.
	\]
	Then,
	\[
	F_{\ga'}(\bZ) ~=~ F_{\ga^*}(\bZ) -\gamma\left(k-\sum_{i\in V}Z_i\right)^2
	\]
	is a nonnegative QPBF by the definition of $\gamma$ and by Lemma \ref{qpbf-xs}. Furthermore, $L_{\ga'}(y^k_Q\mid Q\in\gO^2)\geq 0$ is identical to $L_{\ga^*}(y^k_Q\mid Q\in\gO^2) \geq 0$ by Lemma \ref{identityk2}.
	\qed
	
	\bigskip
	Note that the inequalities introduced by \cite{yang1} correspond to some known extremal nonnegative QPBFs, see \cite{Boros1990} or see Proposition 18 in \cite{Boros2002}.
	Namely, inequality \eqref{Zi} corresponds to $Z_i\geq 0$. Inequality \eqref{ZiZj} corresponds to $Z_iZ_j\geq 0$. Inequality \eqref{ZiZjbar} corresponds to $Z_i(1-Z_j)\geq 0$. Inequality \eqref{ZibarZjbar} corresponds to $(1-Z_i)(1-Z_j)\geq 0$. Inequality \eqref{ZiZjZl+ZibarZjbarZlbar} corresponds to 
	\[
	Z_iZ_jZ_\ell+(1-Z_i)(1-Z_j)(1-Z_\ell) ~=~ 1-Z_i-Z_j-Z_\ell+Z_iZ_j+Z_iZ_\ell+Z_jZ_\ell ~\geq~ 0.
	\]
	Finally, inequality \eqref{ZiZjbarZlbar+ZibarZjZl} corresponds to
	\[
	Z_i(1-Z_j)(1-Z_\ell) + (1-Z_i)Z_jZ_\ell ~=~ Z_i-Z_iZ_j-Z_iZ_\ell + Z_jZ_\ell ~\geq~ 0.
	\]
	These nonnegative quadratic functions are known to belong to $\cG$, see \cite{Boros1992}. There are however many more members in $\cG$ that can be added to tighten up the formulation. In fact the set $\cG$ is partitioned $\cG=\cG^2\cup \cG^3\cup \dots \cup \cG^n$ into families of increasing complexity, where $\cG^d$ denotes the extremal nonnegative quadratic functions that depend on exactly $d$ variables. It is known that none of these classes are empty, and together they contain exponentially many functions \cite{Boros1992}. 
	
	By the above analysis and results of \cite{Boros1992}, the model proposed by \cite{yang1} utilizes exactly the members of $\cG^2\cup \cG^3$. Let us denote by $YAT$ their model and by $LB(YAT)$ and $UB(YAT)$ the corresponding minimum and maximum values.
	
	\begin{equation}\label{YAT-model}
		\begin{array}{crcll}
			\max ~(\min) & 	\displaystyle 	\sum_{Q\in\gO^1} y^1_Q ~+~	\sum_{k=2}^n \frac{1}{\binom{k}{2}}\sum_{Q\in\gO^2} y^k_Q\\
			s.t. \\
			& \displaystyle y^1_{\{i\}} ~+~ \sum_{k=2}^n \frac{1}{k-1} \sum_{Q\in\gO^2\atop i\in Q} y^k_Q &~=~&p_{\{i\}}& \forall ~ i\in V\\
			& \displaystyle \sum_{k=2}^n y^k_Q &~=~&p_Q& \forall ~ Q\in\gO^2\\
			&y^k_{\{i\}}&~\geq~&0&\forall ~ \begin{array}{c}i\in V\\1\leq k\leq n\end{array}\\*[5mm]
			&\displaystyle\sum_{Q\in\Sigma} \ga_Q\cdot y^k_Q &\geq~& 0  &\forall~\begin{array}{c}
				\ga\in\RR^\Sigma \text{ and } \\
				2\leq k\leq n ~\text{where}\\
				F_\ga(\bZ)\in \cG^2\cup \cG^3\end{array} 
		\end{array}
	\end{equation}

	\bigskip
	
	The above theorem provides a surprising connection to the cone of nonnegative quadratic pseudo-Boolean functions, which is known to be equivalent with the cut-polytope (see \cite{DeSimone89}). 
	It was already known that the dual of the maximization problem \eqref{e-Hail} has strong connections to the cut-polytope \cite{Deza92a,Deza92b}, but the above connection to an aggregation of the primal problem is different. 
	
	\bigskip
	Let us denote by $QPB$ the linear programming problem formed by the objective \eqref{obj}, the equalities \eqref{pi} and \eqref{pij}, and the linear inequalites (exponentially many) describing $L_{\cE^2}(\RR^\gO_+)$ as in Theorem \ref{main2}, and denote by $LB(QPB)$ and $UB(QPB)$ the minimum and maximum values of it.

	\begin{equation}\label{A2-model}
		\begin{array}{crcll}
			\max ~(\min) & 	\displaystyle 	\sum_{Q\in\gO^1} y^1_Q ~+~	\sum_{k=2}^n \frac{1}{\binom{k}{2}}\sum_{Q\in\gO^2} y^k_Q\\
			s.t. \\
			& \displaystyle y^1_{\{i\}} ~+~ \sum_{k=2}^n \frac{1}{k-1} \sum_{Q\in\gO^2\atop i\in Q} y^k_Q &~=~&p_{\{i\}}& \forall ~ i\in V\\
			& \displaystyle \sum_{k=2}^n y^k_Q &~=~&p_Q& \forall ~ Q\in\gO^2\\
			&y^k_{\{i\}}&~\geq~&0&\forall ~ \begin{array}{c}i\in V\\1\leq k\leq n\end{array}\\*[5mm]
			&\displaystyle\sum_{Q\in\Sigma} \ga_Q\cdot y^k_Q &\geq~& 0  &\forall~\begin{array}{c}
				\ga\in\RR^\Sigma \text{ and } \\
				2\leq k\leq n ~\text{where}\\
				F_\ga(\bZ)\in \cG\end{array} 
		\end{array}
	\end{equation}
	
	\begin{corollary}\label{c2}
		The aggregation \eqref{A2-model} above is a faithful aggregation of Hailperin's model, that is we have
		\[
		LB(H)~=~LB(QPB) ~~~\text{ and }~~~ UB(H)~=~UB(QPB).
		\]
	\end{corollary}
	\proof
	By Theorem \ref{main2} for every feasible solution $\by$  of problem $QPB$ there exists a corresponding feasible solution $\bx$ of problem \eqref{e-Hail}
	satisfying the equalities \eqref{aijk}.
	\qed
	
	\bigskip
	In practice, we cannot use all inequalities described by Theorem \ref{main2}.
	In this paper we suggest to use some of the functions in $\cG^4\cup\dots\cup \cG^n$ to tighten up the formulation \eqref{obj} - \eqref{ZiZjbarZlbar+ZibarZjZl}. 
	
	For instance, we recall the following result from \cite{Boros1990}: given a subset of the literals $W\subseteq \{Z_1,\dots,Z_n\}\cup\{1-Z_1,\dots,1-Z_n\}$ of cardinality $|W|\geq 4$ and an integer $1\leq \gamma\leq |W|-2$, then the function defined by following binomial expression
	\[
	G_{W,\gamma}(\bZ) ~=~ \binom{\displaystyle\sum_{w\in W}w ~-~ \gamma}{2}
	\]
	is a member of $\cG^{|W|}$. As an example, consider $W=\{Z_1,Z_2,1-Z_3,Z_4\}$ and $\gamma=2$. Then we have

		\begin{align*}
			G_{W,2}(\bZ) 
			&=~ \binom{Z_1+Z_2+(1-Z_3)+Z_4 - 2}{2} \\
			&=~ \frac{(Z_1+Z_2-Z_3+Z_4 - 1)(Z_1+Z_2-Z_3+Z_4 - 2)}{2} \\
			&=~ \frac{(Z_1+Z_2-Z_3+Z_4)^2-3(Z_1+Z_2-Z_3+Z_4)+2}{2}\\
			&=~ 1 -Z_1-Z_2+Z_3-Z_4+Z_1Z_2-Z_1Z_3+Z_1Z_4-Z_2Z_3+Z_2Z_4-Z_3Z_4
		\end{align*}
	Please recall that we have $Z_i^2=Z_i$ and $(1-Z_i)^2=1-Z_i$ for all $i\in V$, since these expressions take only binary values.
	Thus, from this small example we get the inequality
		\[
		y^k_\emptyset 
		-y^k_{\{1\}}-y^k_{\{2\}}+y^k_{\{3\}}-y^k_{\{4\}}
		+y^k_{\{1,2\}}-y^k_{\{1,3\}}+y^k_{\{1,4\}}-y^k_{\{2,3\}}+y^k_{\{2,4\}}-y^k_{\{3,4\}} ~\geq~ 0,
		\]
is valid for $L_\cE(\RR_+^\gO)$.
	Let us introduce such inequalities with $|W|=4$ and $\gamma = 1$. Since we utilize only some functions in $\cG^4$ and none from $\cG^5\cup\dots\cup \cG^n$, we call this weaker model as problem $QPB^{-}$ (and denote by $LB(QPB^{-})$ and $UB(QPB^{-})$ the corresponding minimum and maximum values). Note that $QPB^{-}$ is still a polynomial sized formulation, and thus we can compute the corresponding lower and upper bounds in polynomial time. 
	
	\begin{corollary}\label{c2a2-}
		We have the inequalities
		\[
		\begin{array}{rclll}
			LB(H)&=&LB(QPB)&\geq~ LB(QPB^{-}) &\geq~ LB(YAT) ~~~\text{ and }\\*[3mm]
			UB(H)&=&UB(QPB)&\leq~ UB(QPB^{-}) &\leq~ UB(YAT).
		\end{array}
		\]
		\qed
	\end{corollary}

Since generating all members of $\cG$ is impossible, we demonstrate the improvement by using some functions in $\cG^4$. More specifically, from \cite{Boros1990}, we know for $W\subseteq \{Z_1,\dots,Z_4\}\cup\{1-Z_1,\dots,1-Z_4\}$ the function defined by $G_{W,1}(\bZ)$ is a member of $\cG^4$. In this experiment, we use the row generation methods instead of adding all the inequalities. In the row generation, we first calculate $L_{\ga}(\by^k) = \displaystyle\sum_{Q\in\Sigma} \ga_Q\cdot y^k_Q$ where $F_{\ga}(\bZ)\in \cG^4$. Since $L_{\ga}(\by^k)$ must be nonnegative for any inequality having less than a half of the smallest value $L_{\ga}(\by^k)$ we add the violated row to the model to improve the bounds. We repeat this process until we cannot find any $L_{\ga}(\by^k) \le -0.0001$.

In \cite{Boros14} there is an extensive theoretical comparison of various closed form or polynomially computable lower and upper bound published in the literature and also some new bounds proposed. The $YAT$ and $IPG$ models were introduced later. Accordingly we have as candidates for best lower bounds $LB(YAT)$, $LB(IPG)$, and $LB(Dec)$ introduced in \cite{Boros14}. Similarly, we have $UB(YAT)$, $UB(IPG)$, and $UB(HW)$ as possible best upper bounds.  
In our computational experiments we compare the bounds obtained by our $QPB^-$ model to all these bounds. 

The computational results in Tables \ref{tableQPBUUB} and \ref{tableQPBULB} show that the $QPB^{-}$ model provides significantly improved bounds. 

\begin{table}[h!]
	\small
	\centering
	\begin{tabular}{c c c c}
		\toprule
		\small			n & \tiny$\displaystyle \left( \frac{UB(YAT)-UB(QPB^{-})}{UB(YAT)-UB(H)}\right)$ & \tiny$\displaystyle \left( \frac{UB(IPG)-UB(QPB^{-})}{UB(IPG)-UB(H)}\right)$  & \tiny $\displaystyle \left( \frac{UB(HW)-UB(QPB^{-})}{UB(HW)-UB(H)}\right)$\\
		\midrule
		8	&	60.000 \%	&		99.992\%	&	99.978\%	\\
		9	&	99.658 \%	&		99.973\%	&	99.958\%	\\
		10	&	93.276 \%	&		99.592\%	&	99.207\%	\\
		11	&	84.297 \%	&		98.129\%	&	96.672\%	\\
		12	&	78.174 \%	&		96.498\%	&	94.089\%	\\
		13	&	70.722 \%	&		94.854\%	&	91.264\%	\\
		14	&	69.623 \%	&		94.467\%	&	90.917\%	\\
		15	&	64.098 \%	&		92.467\%	&	87.558\%	\\
		16	&	60.125 \%	&		90.945\%	&	85.467\%	\\
		17	&	56.015 \%	&		89.300\%	&	83.453\%	\\
		18	&	54.434 \%	&		88.816\%	&	82.541\%	\\
		19	&	50.972 \%	&		86.848\%	&	79.690\%	\\
		20	&	50.767 \%	&		86.318\%	&	79.098\%	\\
		21	&	48.104 \%	&		84.998\%	&	77.136\%	\\
		\bottomrule
	\end{tabular}
	\caption{The mean of the percentage relative error improvement over 50 randomly generated problems for each size comparing $UB(QPB^{-})$ to $UB(YAT)$, $UB(IPG)$, and $UB(HW)$.
		\label{tableQPBUUB}}
\end{table}

\begin{table}[h!]
	\small
	\centering
	\begin{tabular}{c c c c}
		\toprule
		\small			n & \tiny$\displaystyle \left( \frac{LB(QPB^{-})-LB(YAT)}{LB(H)-LB(YAT)}\right)$  & \tiny$\displaystyle \left( \frac{LB(QPB^{-})-LB(IPG)}{LB(H)-LB(IPG)}\right)$ & \tiny$\displaystyle \left( \frac{LB(QPB^{-})-LB(Dec)}{LB(H)-LB(Dec)}\right)$\\
		\midrule
		8	&	4.000\%	&		70.000\%	&	70.000\%	\\
		9	&	17.000\%	&		89.727\%	&	89.727\%	\\
		10	&	30.703\%	&		94.500\%	&	94.500\%	\\
		11	&	61.710\%	&		88.405\%	&	88.405\%	\\
		12	&	42.268\%	&		79.562\%	&	79.562\%	\\
		13	&	37.004\%	&		70.801\%	&	70.801\%	\\
		14	&	18.221\%	&		58.238\%	&	58.238\%	\\
		15	&	20.266\%	&		63.113\%	&	63.113\%	\\
		16	&	28.751\%	&		66.856\%	&	66.856\%	\\
		17	&	25.585\%	&		58.996\%	&	58.996\%	\\
		18	&	54.775\%	&		90.172\%	&	90.172\%	\\
		19	&	54.293\%	&		84.425\%	&	84.425\%	\\
		20	&	53.324\%	&		76.771\%	&	76.771\%	\\
		21	&	44.015\%	&		68.799\%	&	68.799\%	\\
		\bottomrule
	\end{tabular}
	\caption{The mean of the percentage relative error improvement over 50 randomly generated problems for each size comparing $LB(QPB^{-})$ to $LB(YAT)$, $LB(IPG)$, and $LB(Dec)$. 
		\label{tableQPBULB}}
\end{table}

	\bigskip
	\bigskip

\section{Conclusions}

In this paper we considered the variable aggregation proposed in \cite{yang1} and observed that the image of the positive orthant under this aggregation is only a convex subcone of the positive orthant of the lower dimensional formulation. Our main result is a complete polyhedral characterization of this subcone that provided a tightening of the aggregated model, both for the minimization and maximization versions, and resulted in significantly improved polynomially computable lower and upper bounds for the probability of the union of events. Our result also make a strong connection between these union bounding models and the cone of the nonnegative quadratic pseudo-Boolean functions. Such a connection was already known for one of the duals of Hailperin's model (see \cite{Deza92a,Deza92b,DezaLau97}). This new connection however seems to be independent of the earlier results, and seems to suggest some fundamental properties of the variables aggregations scheme proposed by \cite{yang1} with respect to the union bounding problem.

We would like to close with a couple of related remarks.

\begin{remark}
	Let us note first that in \cite{yang1}
	the authors conjectured that the \eqref{YAT-model} $YAT$ model provides best possible bounds for $n\leq 7$. While this is true for $n\leq 5$ due to \cite{Grishukin90} and \Cref{main2}, it is not true for $n=6,7$. We include examples, found by random search, for $n=6$ in the appendix. However,
	our computations show that the average relative improvements for $n=6,7$ are quite small, so in our computational tables we included results only for $n\geq 8$.
\end{remark}

\begin{remark}
	Let us remark next that in \cite{YAT16SPL} the authors consider a general scheme where instead of the input parameters $(p_Q\mid Q\in \gO^1\cup\gO^2)$ weighted linear combinations of those are used. This idea is a generalization of the PG model \cite{PG} we cited earlier.
	We would like to add that since this variables aggregation resulted in linear equations for the input parameters, and weighted sum of those is also linearly expressible in terms of the same aggregated variables. Consequently, our inequalities describing the image of the positive orthant can also be added to such models, possibly tightening them up. It may be an interesting future research direction to investigate such weighted tightened models, and their relation to some of the classical bounds from the literature.
\end{remark}

\begin{remark}
	One can notice that the last two columns in Table \ref{tableQPBULB} are identical. The reason is that in all examples we tested (not only in the ones appearing in this table) we found $L(IPG)$ and $LB(Dec)$ agreeing up to six or seven digits. Based on this we conjecture that these two bounds are in fact identical.
\end{remark}

	\bibliographystyle{plain}

	\bibliography{ProbBoundsQPBF}
	
\appendix
\bigskip
\section*{Counter examples to the conjecture of \cite{yang1}:}

For the lower bound case the following input for $n=6$ 
has $LB(YAT)=0.751845$ while $LB(H)=0.758498$ which is significantly different.
\begin{table}[h!]
	\small
	\centering
	\begin{tabular}{c c}
		\toprule
		Set (S) & $P_S$\\
		\midrule
		\{0\}	&	0.313538    	\\
		\{1\}	&	0.31728      	\\
		\{2\}	&	0.269357    	\\
		\{3\}	&	0.32625    	\\
		\{4\}	&	0.315385    	\\
		\{5\}	&	0.291726    	\\
		\{0,1\}	&	0.101524    	\\
		\{0,2\}	&	0.0923267  	\\
		\{0,3\}	&	0.095551    	\\
		\{0,4\}	&	0.103682    	\\
		\{0,5\}	&	0.0874474	\\
		\{1,2\}	&	0.0853875  	\\
		\{1,3\}	&	0.105017    	\\
		\{1,4\}	&	0.106741    	\\
		\{1,5\}	&	0.0974883	\\
		\{2,3\}	&	0.0815002  	\\
		\{2,4\}	&	0.0836789  	\\
		\{2,5\}	&	0.070669	\\
		\{3,4\}	&	0.0955952  	\\
		\{3,5\}	&	0.104005	\\
		\{4,5\}	&	0.10092	\\
		
		\bottomrule
	\end{tabular}
	\caption{Input probabilities for n=6 for lower bound case
		\label{tableInputn6LB}}
\end{table}

For the upper bound case the following input for $n=6$ has 
has $UB(YAT)=0.861996$ while $UB(H)=0.827229$ which is significantly different.
\begin{table}[h!]
	\small
	\centering
	\begin{tabular}{c c}
		\toprule
		Set (S) & $P_S$\\
		\midrule
		\{0\}	&	0.306739    	\\
		\{1\}	&	0.0553987     	\\
		\{2\}	&	0.127482  		\\
		\{3\}	&	0.282613    	\\
		\{4\}	&	0.164909   	\\
		\{5\}	&	0.172431   	\\
		\{0,1\}	&	0.0122853    	\\
		\{0,2\}	&	0.0393019   	\\
		\{0,3\}	&	0.0876492 	\\
		\{0,4\}	&	0.0450683    	\\
		\{0,5\}	&	0.0604373	\\
		\{1,2\}	&	0.00620071   	\\
		\{1,3\}	&	0.0106406  	\\
		\{1,4\}	&	 0.00683587    	\\
		\{1,5\}	&	0.00973235	\\
		\{2,3\}	&	0.0328578  	\\
		\{2,4\}	&	0.0203738 	\\
		\{2,5\}	&	 0.0181882	\\
		\{3,4\}	&	0.0479029  	\\
		\{3,5\}	&	0.0487974	\\
		\{4,5\}	&	0.0214881		\\
		
		\bottomrule
	\end{tabular}
	\caption{Input probabilities for n=6 for upper bound case
		\label{tableInputn6LB}}
\end{table}

\end{document}